\documentclass[12pt]{article}
\usepackage{a4, amssymb, exscale}

\newenvironment{evlist}[2]{
\begin{list}{}{
\setlength{\topsep}{0.5ex plus0.2ex minus0.1ex} 
\setlength{\leftmargin}{#1}
\setlength{\itemsep}{#2 plus0.2ex}
\setlength{\parsep}{0ex plus0.2ex} }}
{\end{list}}

\newcommand{\Realpos}{\mathbb{R}^+}
\newcommand{\Borelpi}{\mathcal{B}^+}
\newcommand{\bMap}{\mathrm{B}}
\newcommand{\uDelta}{\underline{\Delta}}
\newcommand{\Prob}{\mathrm{P}}
\newcommand{\measmap}[2]{#1 {\cdot} #2}
\newcommand{\ext}{\mathrm{ext}\,}
\newcommand{\definition}[1]{\textit{#1}}
\newcommand{\proof}{{\textit{Proof}\enspace}}

\newcommand{\eop}{\ \vbox{\hrule
                       \hbox{\vrule
                             \hskip 6pt
                             \vrule height 6pt width 0pt
                             \vrule}%
                       \hrule}%
                     \vspace{\medskipamount}
                }

\newtheorem{lemma}{Lemma}
\newtheorem{proposition}{Proposition}
\newtheorem{theorem}{Theorem}

\begin{document}

\title{Refining quasi-probability kernels}
\author{Chris Preston}
\date{\small{}}
\maketitle

\begin{quote}
We consider the problem of modifying a quasi-probability kernel in order to improve its properties without 
changing the set of measures whose conditional probabilities it specifies.
\end{quote}

\thispagestyle{empty}

Let $(X,\mathcal{F})$ be some fixed measurable space. We assume throughout that the $\sigma$-algebra $\mathcal{F}$
is countably generated, and emphasise that this property does not automatically carry over to the 
sub-$\sigma$-algebras of $\mathcal{F}$. In fact many sub-$\sigma$-algebras arising naturally in applications 
will fail to be countably generated.

The set of probability measures on $(X,\mathcal{F})$ will be denoted by $\Prob(\mathcal{F})$. A mapping 
$\pi : X \times \mathcal{F} \to \Realpos$ is a \definition{quasi-probability kernel} if $\pi(\,\cdot\,,F)$ 
is an $\mathcal{F}$-measurable mapping for each $F \in \mathcal{F}$ and $\pi(x,\,\cdot\,)$ is a measure on 
$(X,\mathcal{F})$ with $\pi(x,X)$ either $0$ or $1$ for each $x \in X$, and so $\pi(x,\,\cdot\,)$ is either an 
element of $\Prob(\mathcal{F})$ or the zero measure. 

The objects of interest here are quasi-probability kernels, which perhaps raises the question: Why work with 
quasi-probability and not just with probability kernels? One reason is that quasi-probability kernels arise 
naturally, for example  in the study of continuous models or models with unbounded spins in statistical 
mechanics. In such models there is a family of kernels $\{\pi_\Lambda\}$ defined in terms of a Hamiltonian and 
there is no sensible way of defining $\pi_\Lambda(x,\,\cdot\,)$ as a probability measure for all $x \in X$. 
Another reason is that, even when the basic objects of study are probability kernels, many operations result 
in what are really quasi-probability kernels. They are then often artificially modified to turn them into 
probability kernels, but no advantage is incurred by doing this. For the topics to be considered below 
quasi-probability kernels are much easier to deal with than probability kernels, which is another reason
for using them when there are no intrinsic reasons for not doing so.

If $\pi$ is a quasi-probability kernel then the set $\{ x \in X : \pi(x,X) = 1 \}$ is  called the 
\definition{support} of $\pi$ and will be denoted by $S_\pi$. If $\mathcal{E}$ is a sub-$\sigma$-algebra of 
$\mathcal{F}$ then $\pi$ is said to be \definition{$\mathcal{E}$-measurable} if $\pi(\,\cdot\,,F)$ is 
$\mathcal{E}$-measurable for each $F \in \mathcal{F}$. A probability measure $\mu \in \Prob(\mathcal{F})$ is 
said to be \definition{trivial} on $\mathcal{E}$ if $\mu(E) \in \{0,1\}$ for all $E \in \mathcal{E}$.

In what follows let $\mathcal{E}$ be a sub-$\sigma$-algebra of $\mathcal{F}$ and let $\pi$ be an 
$\mathcal{E}$-measurable quasi-probability kernel. Associated with such a kernel $\pi$ is the set 
$\mathcal{J}_{\mathcal{E}}(\pi)$ consisting of those probability measures $\mu$ for which $\pi(\,\cdot\,,F)$ is 
a version the conditional expectation $E_\mu(I_F | \mathcal{E})$ for each $F \in \mathcal{F}$. In other words,
\[\mathcal{J}_{\mathcal{E}}(\pi) = \Bigl\{\, \mu \in \Prob(\mathcal{F}) : \mu(E \cap F) 
= \int I_E\pi(\,\cdot\,,F)\,d\mu \ \mbox{for all}\ E \in \mathcal{E},\, F \in \mathcal{F} \,\Bigr\}\;.\]

In many applications one the main tasks is to analyse the set $\mathcal{J}_\mathcal{E}(\pi)$. For this it is 
useful to be able to exploit additional properties of the kernel $\pi$, and let us first mention a property
that always holds (see Proposition~\ref{prop_wip_e21}): 
\begin{evlist}{24pt}{6pt}
\item[(0)]
If $\mu \in \mathcal{J}_{\mathcal{E}}(\pi)$ is trivial on $\mathcal{E}$ then
$\mu = \pi(x,\,\cdot\,)$ for some $x \in S_\pi$.
\end{evlist}
The most important properties (which do not always hold) are probably those occurring in the following 
definitions:

\begin{evlist}{24pt}{6pt}
\item[(1)]
$\pi$ is said to be \definition{proper} if $\pi(\,\cdot\,,E\cap F) = I_E\pi(\,\cdot\,,F)$ for all 
$E \in \mathcal{E}$, $F \in \mathcal{F}$.

\item[(2)]
$\pi$ will be called \definition{adapted} if $\pi(x,\,\cdot\,) \in \mathcal{J}_{\mathcal{E}}(\pi)$ for all 
$x \in S_\pi$.

\item[(3)]
$\pi$ is called \definition{normal} if it is adapted and $\pi(x,\,\cdot\,)$ is trivial on $\mathcal{E}$ for 
all $x \in S_\pi$.
\end{evlist}

The term \definition{normal} is taken from Dynkin \cite{dynkin1}. Being \definition{adapted} corresponds to what 
Dynkin \cite{dynkin1}, \cite{dynkin2} calls a \definition{$(\mathcal{J}_{\mathcal{E}}(\pi),\mathcal{E})$-kernel}.
It is well-known -- and will be shown in Lemma~\ref{lemma_wip_21} -- that a proper kernel is normal.
If $\pi$ is normal then (0) implies that $\mu \in \mathcal{J}_{\mathcal{E}}(\pi)$
is trivial on $\mathcal{E}$ if and only if $\mu = \pi(x,\,\cdot\,)$ for some $x \in S_\pi$.

The reason why normal kernels are important is because of the following:
Denote the extreme points of the convex set $\mathcal{J}_\mathcal{E}(\pi)$ by $\ext \mathcal{J}_\mathcal{E}(\pi)$.
Then it is well-known (and a proof is provided in Proposition~\ref{prop_wip_31}) that an element of 
$\mathcal{J}_\mathcal{E}(\pi)$ is extreme if and only if it is trivial on $\mathcal{E}$. It follows that if the 
kernel $\pi$ is normal then $\ext \mathcal{J}_{\mathcal{E}}(\pi)$ consists of exactly the elements in
$\mathcal{J}_{\mathcal{E}}(\pi)$ having the form $\pi(x,\,\cdot\,)$ for some $x \in S_\pi$, and it is this fact 
which plays a crucial role in applications. 

Now it may be that $\pi$ itself fails to have one of these three properties but that it is possible to modify 
$\pi$ to obtain a normal (resp.\ adapted resp.\ proper) kernel $\varrho$ such that 
$\mathcal{J}_{\mathcal{E}}(\varrho) = \mathcal{J}_{\mathcal{E}}(\pi)$ holds. The analysis of the set 
$\mathcal{J}_{\mathcal{E}}(\pi)$ can then be carried out using $\varrho$ instead of $\pi$. This is the topic to 
be discussed here, and so let us start by describing how the kernels will be modified.

If $D \in \mathcal{E}$ then putting $\varrho(x,F) = I_D(x)\pi(x,F)$ for all $F \in \mathcal{F}$, $x \in X$
defines an $\mathcal{E}$-measurable quasi-probability kernel $\varrho$ with $S_\varrho = D \cap S_\pi$, which 
will be called the \definition{restriction of $\pi$ to $D$}. An $\mathcal{E}$-measurable quasi-probability 
kernel $\varrho$ will be called a \definition{refinement} of $\pi$ if $\varrho$ is the restriction of $\pi$ to 
$D$ for some $D \in \mathcal{E}$ and $\mathcal{J}_{\mathcal{E}}(\varrho) = \mathcal{J}_{\mathcal{E}}(\pi)$.

If $\mathcal{E}$ is countably generated then it is straightforward to show that there exists a proper refinement 
of $\pi$, and this will done in the proof of Theorem~\ref{theorem_wip_11}. Very similar results for
probability -- rather than quasi-probability -- kernels can be found in Halmos \cite{halmos} and 
Sokal \cite{sokal}. If $\mathcal{E}$ is not countably generated then in general a proper refinement is too much 
to expect and so our aim is to find conditions on $\pi$ which ensure that a normal refinement $\varrho$ of 
$\pi$ exists. An important role  in this endeavour will be played by the set
\[ \mathcal{J}_\star(\pi) = \Bigl\{\, \mu \in \Prob(\mathcal{F}) :  
                    \mu(F) = \int \pi(\,\cdot\,,F)\,d\mu\ \mbox{for all $F \in \mathcal{F}$}\,\Bigr\}\;.\]
Thus $\mathcal{J}_{\mathcal{E}}(\pi)\subset \mathcal{J}_\star(\pi)$, since if 
$\mu \in \mathcal{J}_{\mathcal{E}}(\pi)$ then for all $F \in \mathcal{F}$
\[ \mu(F) = \mu(X \cap F) = \int I_X\pi(\,\cdot\,,F)\,d\mu = \int \pi(\,\cdot\,,F)\,d\mu\;,\] 
with equality when $\pi$ is proper, since in this case if $\mu \in \mathcal{J}_\star(\pi)$ then
\[ \mu(E \cap F) = \int \pi(\,\cdot\,,E \cap F)\,d\mu = \int I_E \pi(\,\cdot\,,F)\,d\mu\]
for all $E \in \mathcal{E}$, $F \in \mathcal{F}$. The main result (Theorem~\ref{theorem_wip_21}) states that if 
$\mathcal{J}_{\mathcal{E}}(\pi) = \mathcal{J}_\star(\pi)$ then there does exist a normal refinement of $\pi$. 

The proof of Theorem~\ref{theorem_wip_21} makes use of a fact first noted by Blackwell and Dubins in 
\cite{blackwelldubins}: There is a least sub-$\sigma$-algebra $\mathcal{S}_\pi$ such that $\pi$ is 
$\mathcal{S}_\pi$-measurable (namely the $\sigma$-algebra generated by the mappings $\pi(\,\cdot\,,F)$, 
$F \in \mathcal{F}$) and the assumption that $\mathcal{F}$ is countably generated implies that $\mathcal{S}_\pi$ 
is also countably generated. We can thus consider $\pi$ as an $\mathcal{S}_\pi$-measurable kernel and apply 
Theorem~\ref{theorem_wip_11} to obtain a proper and hence normal refinement $\varrho$ of $\pi$. Of course, at 
first glance $\varrho$ is only normal as an $\mathcal{S}_\pi$-measurable kernel, but 
$\mathcal{J}_{\mathcal{E}}(\pi) \subset \mathcal{J}_{\mathcal{S}_\pi}(\pi)$, since
$\mathcal{S}_\pi \subset \mathcal{E}$, and thus 
$\mathcal{J}_{\mathcal{E}}(\pi) = \mathcal{J}_{\mathcal{S}_\pi}(\pi)$, since 
$\mathcal{J}_{\mathcal{S}_\pi}(\pi) \subset \mathcal{J}_\star(\pi)$ and 
$\mathcal{J}_{\mathcal{E}}(\pi) = \mathcal{J}_\star(\pi)$. From this it will follow that $\varrho$ is also normal 
as an $\mathcal{E}$-measurable kernel.

Before going any further let us introduce some more convenient notation for measures and kernels. For each 
sub-$\sigma$-algebra $\mathcal{F}_0$ of $\mathcal{F}$ denote the set of bounded $\mathcal{F}_0$-measurable 
mappings from $X$ to $\Realpos$ by $\bMap(\mathcal{F}_0)$. If $\mu$ is a finite measure on $(X,\mathcal{F})$ then 
for each $f \in \bMap(\mathcal{F})$ we write $\mu(f)$ instead of $\int f \,d\mu$. The measure is thus considered 
as a mapping $\mu : \bMap(\mathcal{F}) \to \Realpos$, and what was previously $\mu(F)$ now becomes $\mu(I_F)$.

In the same way, let $\varrho : X \times \mathcal{F} \to \Realpos$ be a bounded kernel, meaning that 
$\varrho(\,\cdot\,,F)$ is an $\mathcal{F}$-measurable mapping for each $F \in \mathcal{F}$ and there exists
$N \ge 0$ such that $\varrho(x,\,\cdot\,)$ is a finite measure on $(X,\mathcal{F})$ with $\varrho(x,X) \le N$ 
for each $x \in X$. Then for each $f \in \bMap(\mathcal{F})$ and each $x \in X$ we write $\varrho(f)(x)$ instead
of $\int f(y)\,\varrho(x,dy)$. Thus, since the mapping $x \mapsto \int f(y) \varrho(x,dy)$ defines an element of 
$\bMap(\mathcal{F})$, the kernel is considered as a mapping $\varrho : \bMap(\mathcal{F}) \to \bMap(\mathcal{F})$,
and what was previously $\varrho(x,F)$ now becomes $\varrho(I_F)(x)$.

Let $\tau,\,\varrho : \bMap(\mathcal{F}) \to \bMap(\mathcal{F})$ be bounded kernels; then the bounded kernel 
$\varrho\tau$, which using the old notation is defined by $(\varrho\tau)(x,F) = \int \tau(y,F) \varrho(x,dy)$, 
is now given as a mapping $\varrho\tau : \bMap(\mathcal{F}) \to \bMap(\mathcal{F})$ by
$(\varrho\tau)(f) = \varrho(\tau(f))$, and so the `product' of the kernels is just functional composition. 
Moreover, if $\mu$ is a finite  measure then the finite measure $\mu\tau : \bMap(\mathcal{F}) \to \Realpos$,
which using the old notation is defined by $(\mu\tau)(F) = \int \tau(\,\cdot\,,F)\,d\mu$, is now given as a 
mapping $\mu\tau : \bMap(\mathcal{F}) \to \Realpos$ by $(\mu\tau)(f) = \mu(\tau(f))$, which is again functional 
composition. In particular, there are  no problems with the associativity of the various operations, since 
this holds trivially for the composition of mappings.

The $\mathcal{E}$-measurable quasi-probability kernel $\pi$ will thus now be considered as a mapping
$\pi : \bMap(\mathcal{F}) \to \bMap(\mathcal{F})$. It is easily checked that $\pi$ being $\mathcal{E}$-measurable
is the same as having $\pi(f) \in \bMap(\mathcal{E})$ for all $f \in \bMap(\mathcal{F})$, and that $\pi$ will be 
proper if and only if $\pi(gf) = g\pi(f)$ for all $g \in \bMap(\mathcal{E})$, $f \in \bMap(\mathcal{F})$.
Moreover,
\[\mathcal{J}_{\mathcal{E}}(\pi) = \{ \mu \in \Prob(\mathcal{F}) : \mu(gf) = \mu(g\pi(f))
\ \mbox{for all}\ g \in \bMap(\mathcal{E}),\, f \in \bMap(\mathcal{F}) \}\;,\]
$\mathcal{J}_\star(\pi) = \{ \mu \in \Prob(\mathcal{F}) : \mu\pi = \mu \}$ and the probability measure $\mu$ is 
trivial on $\mathcal{E}$ if and only if $\mu(gf) = \mu(g)\mu(f)$ for all $g \in \bMap(\mathcal{E})$, 
$f \in \bMap(\mathcal{F})$.

For each $x \in X$ let $\varepsilon_x \in \Prob(\mathcal{F})$ be the point mass at $x$, and thus
$\varepsilon_x(f) = f(x)$ for all $f \in \bMap(\mathcal{F})$. For each $x \in X$ there is the measure 
$\varepsilon_x\pi$, and by definition $(\varepsilon_x\pi)(f) = \pi(f)(x)$ for all $f \in \bMap(\mathcal{F})$. 
With the previous notation this means that $(\varepsilon_x\pi)(I_F) = \pi(x,F)$. In particular, $\pi$ will be 
normal if and only if $\varepsilon_x\pi \in \mathcal{J}_{\mathcal{E}}(\pi)$ with $\varepsilon_x\pi$ trivial on 
$\mathcal{E}$ for all $x \in S_\pi$.

We will often make use of the fact that $\mu(I_{S_\pi}) = 1$ for all $\mu \in \mathcal{J}_{\mathcal{E}}(\pi)$, which 
holds since $\mu(I_{S_\pi}) = \mu(\pi(1)) = (\mu\pi)(1) = \mu(1) = 1$. Note that if $D \in \mathcal{E}$ then the 
restriction $\varrho$ of $\pi$ to $D$ is given by $\varrho(f) = I_D\pi(f)$ for all $f \in \bMap(\mathcal{F})$.

\begin{lemma}\label{lemma_wip_11}
Let $D \in \mathcal{E}$ and $\varrho$ be the restriction of $\pi$ to $D$. Then $\varrho$ is a refinement of $\pi$ 
if and only if $\mu(I_D) = 1$ for all $\mu \in \mathcal{J}_{\mathcal{E}}(\pi)$. 
\end{lemma}

\proof 
Let us suppose first that $\mu(I_D) = 1$ for all $\mu \in \mathcal{J}_{\mathcal{E}}(\pi)$. If 
$\mu \in \mathcal{J}_{\mathcal{E}}(\pi)$ then $\mu(g\varrho(f)) = \mu(gI_D\pi(f)) = \mu(g\pi(f)) = \mu(gf)$ for 
all $f \in \bMap(\mathcal{F})$, $g \in \bMap(\mathcal{E})$, since $\mu(I_D) = 1$, and therefore
$\mu \in \mathcal{J}_{\mathcal{E}}(\varrho)$. On the other hand, if $\mu \in \mathcal{J}_{\mathcal{E}}(\varrho)$ 
then $\mu(I_D) = 1$ (since $D \supset S_\varrho$ and $\mu(I_{S_\varrho}) = 1$) and so for all 
$f \in \bMap(\mathcal{F})$, $g \in \bMap(\mathcal{E})$ we have 
$\mu(g\pi(f)) = \mu(gI_D\pi(f)) = \mu(g\varrho(f)) = \mu(gf)$. 
Hence $\mu \in \mathcal{J}_{\mathcal{E}}(\pi)$ and thus 
$\mathcal{J}_{\mathcal{E}}(\varrho) = \mathcal{J}_{\mathcal{E}}(\pi)$, which shows that $\varrho$ is a refinement 
of $\pi$. Suppose conversely that $\varrho$ is a refinement of $\pi$; then, since $D  \supset S_\varrho$, it 
follows that $\mu(I_D) = 1$ for all $\mu \in \mathcal{J}_{\mathcal{E}}(\varrho) = \mathcal{J}_{\mathcal{E}}(\pi)$. 
\eop

\begin{lemma}\label{lemma_wip_21}
If $\pi$ is proper then it is also normal.
\end{lemma}

\proof
Let $x \in S_\pi$; then $\pi(1)(x) = 1$ and hence for all $f \in \bMap(\mathcal{F})$, $g \in \bMap(\mathcal{E})$
\begin{eqnarray*}
(\varepsilon_x\pi)(g\pi(f)) &=& \pi(g\pi(f))(x)\\ &=& 
g(x)\pi(f)(x)\pi(1)(x) = g(x)\pi(f)(x) = \pi(gf)(x) = (\varepsilon_x\pi)(gf)
\end{eqnarray*}
Thus $\varepsilon_x\pi \in \mathcal{J}_{\mathcal{E}}(\pi)$. Moreover,
$(\varepsilon_x\pi)(I_E) = \pi(I_E)(x) = I_E(x)\pi(1)(x) \in \{0,1\}$ for all $E \in \mathcal{E}$. \eop

The next proposition gives some properties of adapted kernels.

\begin{proposition}\label{prop_wip_11}
Suppose $\pi$ is adapted; then:
\begin{evlist}{10pt}{6pt}
\item[(1)]
$\mathcal{J}_{\mathcal{E}}(\pi) = \mathcal{J}_\star(\pi)$.

\item[(2)]
If $\tau$ is an arbitrary (not necessarily $\mathcal{E}$-measurable) quasi-probability kernel
\hphantom{xxx}such that $\mathcal{J}_{\mathcal{E}}(\pi) \subset \mathcal{J}_\star(\tau)$ then $\pi\tau = \pi$.
In particular, $\pi\pi = \pi$. 
\end{evlist}
\end{proposition}

\proof 
(1)\enskip 
Since $\pi$ is adapted it follows that $\pi(g\pi(f)) = \pi(gf)$ for all $f \in \bMap(\mathcal{F})$, 
$g \in \bMap(\mathcal{E})$, since if $x \notin S_\pi$ then $\pi(g\pi(f))(x) = 0 = \pi(gf)(x)$. Thus if 
$\mu \in \mathcal{J}_\star(\pi)$ then for all $f \in \bMap(\mathcal{F})$, $g \in \bMap(\mathcal{E})$
\[ \mu(g\pi(f)) = (\mu\pi)(g\pi(f)) = \mu(\pi(g\pi(f))) = \mu(\pi(gf)) = (\mu\pi)(gf) = \mu(gf) \]
and so $\mu \in \mathcal{J}_{\mathcal{E}}(\pi)$. This shows that 
$\mathcal{J}_\star(\pi) \subset \mathcal{J}_{\mathcal{E}}(\pi)$ and hence that
$\mathcal{J}_{\mathcal{E}}(\pi) = \mathcal{J}_\star(\pi)$.

(2)\enskip 
Let $\tau$ be an arbitrary quasi-probability kernel with 
$\mathcal{J}_{\mathcal{E}}(\pi) \subset \mathcal{J}_\star(\tau)$. If $x \in S_\pi$ then 
$\varepsilon_x\pi \in \mathcal{J}_{\mathcal{E}}(\pi) \subset \mathcal{J}_\star(\tau)$ and so for all 
$f \in \bMap(\mathcal{F})$
\[ (\pi\tau)(f)(x) = 
(\varepsilon_x(\pi\tau))(f) = ((\varepsilon_x\pi)\tau)(f) =  (\varepsilon_x\pi)(f) = \pi(f)(x)\;.\]
On the other hand, if $x \notin S_\pi$ then $\varepsilon_x\pi = 0$, hence
\[ (\pi\tau)(f)(x) = 
(\varepsilon_x(\pi\tau))(f) = ((\varepsilon_x\pi)\tau)(f) =  0 = (\varepsilon_x\pi)(f) = \pi(f)(x)\]
and therefore $(\pi\tau)(f)(x) = \pi(f)(x)$ for all $f \in \bMap(\mathcal{F})$, $x \in X$, which means that
$\pi\tau = \pi$. In particular $\pi\pi = \pi$, since 
$\mathcal{J}_{\mathcal{E}}(\pi) \subset \mathcal{J}_\star(\pi)$.
\eop

Note that if $\varrho$ is an adapted refinement of $\pi$ then by Proposition~\ref{prop_wip_11}
$\varrho\pi = \varrho$, since here
$\mathcal{J}_{\mathcal{E}}(\varrho) = \mathcal{J}_{\mathcal{E}}(\pi) \subset \mathcal{J}_\star(\pi)$.

Since $\mathcal{F}$ is countably generated there exists a countable subset $G$ of $\bMap(\mathcal{F})$ which 
determines finite measures in that if $\mu_1,\,\mu_2$ are finite measures on $(X,\mathcal{F})$ then 
$\mu_1 = \mu_2$ if and only if $\mu_1(f) = \mu_2(f)$ for all $f \in G$. For example, there exists a countable 
algebra $\mathcal{A}$ with $\sigma(\mathcal{A}) = \mathcal{F}$ and then $\{ I_A : A \in \mathcal{A} \}$ has 
this property. In what follows let $G$ be such a countable determining set. 

We now come to the first result about the existence of refinements.

\begin{theorem}\label{theorem_wip_11}
If the sub-$\sigma$-algebra $\mathcal{E}$ is countably generated then there exists a proper 
$\mathcal{E}$-measurable refinement $\varrho$ of $\pi$.
\end{theorem}

\proof 
As well as the countable determining set $G$ choose a countable determining set $G'$ for finite measures on 
$(X,\mathcal{E})$. Let
\[ D = \{ x \in S_\pi : \pi(gf)(x) = g(x)\pi(f)(x)\ \mbox{for all $g \in G'$,  $f \in G$} \}\;;\]
then $D = \bigcap_{f \in G} \bigcap_{g \in G'} D_{f,g}$, where 
$D_{f,g} =  \{ x \in S_\pi: \pi(gf)(x) = g(x)\pi(f)(x) \}$. In particular $D \in \mathcal{E}$. Let 
$\mu \in \mathcal{J}_{\mathcal{E}}(\pi)$, $g \in G'$, $f \in G$; then for all $h \in \bMap(\mathcal{E})$ we have
$\mu(h\pi(gf)) = \mu(hgf)) = \mu(hg\pi(f))$, and thus $\mu(I_{D_{g,h}}) = 1$ since $\mu(I_{S_\pi}) = 1$. Hence 
$\mu(I_D) = 1$, since $G \times G'$ is countable. 

Now let $\varrho$ be the restriction of $\pi$ to $D$; then by Lemma~\ref{lemma_wip_11} $\varrho$ is a refinement 
of $\pi$ since $\mu(I_D) = 1$ for each $\mu \in \mathcal{J}_{\mathcal{E}}(\pi)$. Fix $x \in D$ and $f \in G$; then
\[(\varepsilon_x\varrho)(gf) = \varrho(gf)(x) = \pi(gf)(x) = g(x)\pi(f)(x) = g(x)\varrho(f)(x) 
= g(x)(\varepsilon_x\varrho)(f)\;,\]
i.e., $(\varepsilon_x\varrho)(gf) = g(x)(\varepsilon_x\varrho)(f)$ for all $g \in G'$ and hence for all 
$g \in \bMap(\mathcal{E})$, since $g \mapsto (\varepsilon_x\varrho)(gf)$ and 
$g \mapsto g(x)(\varepsilon_x\varrho)(f)$ are both finite measures on $(X,\mathcal{E})$. From this it follows 
that $(\varepsilon_x\varrho)(gf) = g(x)(\varepsilon_x\varrho)(f)$ for all $g \in \bMap(\mathcal{E})$, 
$f \in \bMap(\mathcal{F})$, since $f \mapsto (\varepsilon_x\varrho)(gf)$ and 
$f \mapsto g(x)(\varepsilon_x\varrho)(f)$ are both finite measures on $(X,\mathcal{F})$. This shows that 
$\varrho(gf)(x) = g(x)\varrho(f)(x)$ for all $g \in \bMap(\mathcal{E})$, $f \in \bMap(\mathcal{F})$, $x \in D$.
But if $x \notin D$ then $\varrho(gf)(x) = 0 = g(x)\varrho(f)(x)$ and therefore $\varrho$ is a proper 
$\mathcal{E}$-measurable kernel. \eop

For each $\mu \in \Prob(\mathcal{F})$ put $\Delta_\mu^\pi = \{ x \in S_\pi : \varepsilon_x\pi = \mu \}$. One 
reason for requiring $\mathcal{F}$ to be countably generated is that it ensures the measurability of sets such 
as $\Delta_\mu^\pi$:

\begin{lemma}\label{lemma_wip_31}
Let $\mu \in \Prob(\mathcal{F})$; then
$\Delta_\mu^\pi  = \bigcap_{f\in G} \{ x \in S_\pi : \pi(f)(x) = \mu(f) \}$
and so in particular $\Delta_\mu^\pi \in \mathcal{E}$.
\end{lemma}

\proof 
Let $x \in S_\pi$; then
\[
\Delta_\mu^\pi
 = \{ x \in S_\pi : \mbox{$(\varepsilon_x\pi)(f) = \mu(f)$ for all $f \in G$} \}
 = \bigcap_{f\in G} \{ x \in S_\pi : \pi(f)(x) = \mu(f) \}
\]
and $\{ x \in S_\pi : \pi(f)(x) = \mu(f) \} \in \mathcal{E}$ for each $f \in G$.
\eop

To increase the legibility we use $\uDelta_\mu^\pi$ to denote the indicator function $I_{\Delta_\mu^\pi}$ of 
$\Delta_\mu^\pi$; by Lemma~\ref{lemma_wip_31} $\uDelta_\mu^\pi \in \bMap(\mathcal{E})$. The following is taken 
from Dynkin \cite{dynkin1}:

\begin{lemma}\label{lemma_wip_41}
If $\mu \in \mathcal{J}_{\mathcal{E}}(\pi)$ then $\mu(\uDelta^\pi_\mu) = 1$ if and only if $\mu$ is trivial on 
$\mathcal{E}$.
\end{lemma}

\proof
For each $f \in \bMap(\mathcal{F})$ let $D_f = \{ x \in X : \pi(f)(x) = \mu(f) \}$. If $f \in \bMap(\mathcal{F})$
and $\mu \in \Prob(\mathcal{F})$ then $\mu(I_{D_f}) = 1$ if and only if 
$\mu(g\pi(f)) = \mu(g\mu(f)) = \mu(g)\mu(f)$ for all $g \in \bMap(\mathcal{E})$, since 
$\pi(f) \in \bMap(\mathcal{E})$. Now let $\mu \in \mathcal{J}_{\mathcal{E}}(\pi)$; then $\mu(I_{D_f}) = 1$ if and 
only if $\mu(gf) = \mu(g)\mu(f)$ for all $g \in \bMap(\mathcal{E})$, since here $\mu(g\pi(f)) = \mu(gf)$.

Suppose $\mu$ is trivial on $\mathcal{E}$; then $\mu(gf) = \mu(g)\mu(f)$ does hold for all 
$g \in \bMap(\mathcal{E})$, and hence $\mu(I_{D_f}) = 1$ for all $f \in \bMap(\mathcal{F})$. But by 
Lemma~\ref{lemma_wip_31} $\Delta_\mu^\pi = \bigcap_{f \in G} D_f$ and hence $\mu(\uDelta_\mu^\pi) = 1$, since 
$G$ is countable.

Suppose conversely that $\mu(\uDelta_\mu^\pi) = 1$, let $E \in \mathcal{E}$ and put $h = I_E$; then 
$\Delta_\mu^\pi \subset D_h$ and so $\mu(I_{D_h}) = 1$, which implies that $\mu(gh) = \mu(g)\mu(h)$ for all 
$g \in \bMap(\mathcal{E})$. In particular, with $g = h$, it follows that $\mu(h^2) = (\mu(h))^2$, i.e., 
$\mu(I_E) = (\mu(I_E))^2$ and therefore $\mu(I_E) \in \{0,1\}$. This shows that $\mu$ is trivial on 
$\mathcal{E}$. \eop

In the proof of Lemma~\ref{lemma_wip_41} we used the fact that if $\mathcal{F}_0$ is a sub-$\sigma$-algebra
of $\mathcal{F}$ and $f_1,\,f_2 \in \bMap(\mathcal{F}_0)$, $\mu \in \Prob(\mathcal{F})$ with
$\mu(gf_1) = \mu(gf_2)$ for all $g \in \bMap(\mathcal{F}_0)$ then $\mu(I_D) = 1$, where 
$D = \{ x \in X : f_1(x) = f_2(x) \}$. This will also be made use of several times below.

\begin{proposition}\label{prop_wip_e21}
If $\mu \in \mathcal{J}_{\mathcal{E}}(\pi)$ is trivial on $\mathcal{E}$ then $\mu = \varepsilon_x\pi$ for some 
$x \in S_\pi$. In particular, if $\pi$ is normal then $\mu \in \mathcal{J}_{\mathcal{E}}(\pi)$ is trivial on 
$\mathcal{E}$ if and only if $\mu = \varepsilon_x\pi$ for some $x \in S_\pi$.
\end{proposition}

\proof 
This follows immediately from Lemma~\ref{lemma_wip_41}.
\eop

For each $x \in S_\pi$ let $\Delta^\pi_x = \{ y \in S_\pi : \varepsilon_y\pi = \varepsilon_x\pi \}$; thus
$\Delta^\pi_x = \Delta^\pi_{\varepsilon_x\pi}$. Also denote the indicator function $I_{\Delta_x^\pi}$ by 
$\uDelta_x^\pi$.

\begin{proposition}\label{prop_wip_21}
$\pi$ is normal if and only if it is adapted and $(\varepsilon_x\pi)(\uDelta_x^\pi) = 1$ for all $x \in S_\pi$.
\end{proposition}

\proof 
This follows immediately from Lemma~\ref{lemma_wip_41}.
\eop

\begin{theorem}\label{theorem_wip_21}
If $\mathcal{J}_{\mathcal{E}}(\pi) = \mathcal{J}_\star(\pi)$ then there exists a normal refinement $\varrho$ of 
$\pi$.
\end{theorem}

\proof
There is a least sub-$\sigma$-algebra $\mathcal{S}_\pi$ of $\mathcal{F}$ such that $\pi$ is 
$\mathcal{S}_\pi$-measurable, this being the intersection of all such sub-$\sigma$-algebras; thus 
$\mathcal{S}_\pi \subset \mathcal{E}$. The next fact (taken from Theorem 1 in Blackwell and 
Dubins \cite{blackwelldubins}) follows from the assumption that $\mathcal{F}$ is countably generated.

\begin{lemma}\label{lemma_wip_51}
The sub-$\sigma$-algebra $\mathcal{S}_\pi$ is countably generated. 
\end{lemma}

\proof 
Since $\mathcal{F}$ is countably generated there exists a countable algebra $\mathcal{A}$ with
$\mathcal{F} = \sigma(\mathcal{A})$. Let $\mathcal{S}'_\pi$ be the least sub-$\sigma$-algebra of $\mathcal{F}$ 
such that $\pi(I_A) \in \mathcal{S}'_\pi$ for all $A \in \mathcal{A}$. Then 
$\mathcal{S}'_\pi \subset \mathcal{S}_\pi$ and  $\mathcal{S}'_\pi$ is countably generated, since 
$\{ \pi(I_A) : A \in \mathcal{A} \}$ is a countable set of mappings. Let 
$\mathcal{F}' = \{ F \in \mathcal{F} : \pi(I_F) \in \mathcal{S}'_\pi \}$; then $\mathcal{F}'$ contains 
$\mathcal{A}$ and is a monotone class and hence by the monotone class theorem $\mathcal{F}' = \mathcal{F}$. 
Thus $\pi$ is $\mathcal{S}'_\pi$-measurable and so $\mathcal{S}_\pi \subset \mathcal{S}'_\pi$, i.e., 
$\mathcal{S}_\pi = \mathcal{S}'_\pi$, which shows that $\mathcal{S}_\pi$ is countably generated.
\eop

Since $\mathcal{S}_\pi \subset \mathcal{E}$ it follows immediately that
$\mathcal{J}_{\mathcal{E}}(\pi) \subset \mathcal{J}_{\mathcal{S}_\pi}(\pi)$. Together with the assumption that 
$\mathcal{J}_{\mathcal{E}}(\pi) = \mathcal{J}_\star(\pi)$ this gives us
$\mathcal{J}_{\mathcal{E}}(\pi) = \mathcal{J}_{\mathcal{S}_\pi}(\pi)$. 

Consider $\pi$ as an $\mathcal{S}_\pi$-measurable kernel; then by Lemma~\ref{lemma_wip_51} and
Theorem~\ref{theorem_wip_11} there exists a proper refinement $\varrho$ of $\pi$. More precisely: There exists 
$D \in \mathcal{S}_\pi$ such that $\varrho$ is the restriction of $\pi$ to $D$,
$\mathcal{J}_{\mathcal{S}_\pi}(\varrho) = \mathcal{J}_{\mathcal{S}_\pi}(\pi)$ and such that $\varrho$ is proper as 
an $\mathcal{S}_\pi$-measurable kernel. By Lemma~\ref{lemma_wip_11} $\mu(I_D) = 1$ for all 
$\mu \in \mathcal{J}_{\mathcal{S}_\pi}(\pi)$ and by Lemma~\ref{lemma_wip_21} $\varrho$ is a normal 
$\mathcal{S}_\pi$-measurable kernel.

Now consider $\pi$ and $\varrho$ as $\mathcal{E}$-measurable kernels. Then by Lemma~\ref{lemma_wip_11} $\varrho$
is still a refinement of $\pi$, since $D \in \mathcal{E}$ and 
$\mathcal{J}_{\mathcal{E}}(\pi) = \mathcal{J}_{\mathcal{S}_\pi}(\pi)$. Therefore 
$\mathcal{J}_{\mathcal{E}}(\varrho) = \mathcal{J}_{\mathcal{E}}(\pi)$, and so
$\mathcal{J}_{\mathcal{E}}(\varrho) = \mathcal{J}_{\mathcal{E}}(\pi) = \mathcal{J}_{\mathcal{S}_\pi}(\pi)
 = \mathcal{J}_{\mathcal{S}_\pi}(\varrho)$,
i.e., $\mathcal{J}_{\mathcal{E}}(\varrho) = \mathcal{J}_{\mathcal{S}_\pi}(\varrho)$. It follows that $\varrho$ is 
an adapted  $\mathcal{E}$-measurable kernel and hence by Proposition~\ref{prop_wip_21} $\varrho$ is a normal  
$\mathcal{E}$-measurable kernel, since the condition $(\varepsilon_x\varrho)(\uDelta_x^\varrho) = 1$ for all 
$x \in S_\varrho$ does not depend on which of the sub-$\sigma$-algebras $\mathcal{S}_\pi$ and $\mathcal{E}$ is 
being used.

This completes the proof of Theorem~\ref{theorem_wip_21}. \eop

Below we will need the following simple fact about proper kernels:

\begin{lemma}\label{lemma_wip_61}
$\pi$ is proper if and only if $\pi(I_E) = I_E\pi(1)$ for all $E \in \mathcal{E}$.
\end{lemma}

\proof 
If $\pi$ is proper then $\pi(I_E) = \pi(I_E1) = I_E\pi(1)$ for all $E \in \mathcal{E}$. Thus suppose conversely 
that $\pi(I_E) = I_E\pi(1)$ for all $E \in \mathcal{E}$. Let $F \in \mathcal{F}$ and $E \in \mathcal{E}$; then
$\pi(I_EI_F) \le \pi(I_E) = I_E\pi(1)$ and $\pi(I_EI_F) \le \pi(I_F)$, since $I_EI_F \le \min\{I_E,I_F\}$;
thus $\pi(I_EI_F) \le I_E\pi(1)\pi(I_F) = I_E\pi(I_F)$. In the same way
$\pi(I_{X \setminus E}I_F) \le I_{X \setminus E}\pi(I_F)$. But
$\pi(I_E I_F) + \pi(I_{X \setminus E}I_F) = \pi(I_F) + I_E\pi(I_F) + I_{X \setminus E}\pi(I_F)$
and so in particular $\pi(I_EI_F) = I_E\pi(I_F)$. Hence $\pi$ is proper. \eop

We now consider some conditions which are equivalent to the kernel $\pi$ being normal. As in the proof of 
Theorem~\ref{theorem_wip_21} let $\mathcal{S}_\pi$  denote the least sub-$\sigma$-algebra of $\mathcal{F}$ 
such that $\pi$ is $\mathcal{S}_\pi$-measurable. Since $\pi$ is $\mathcal{S}_\pi$-measurable, it follows from
Lemma~\ref{lemma_wip_31} that $\Delta^\pi_x \in \mathcal{S}_\pi$. Note that $x \in \Delta_x^\pi$ for all 
$x \in S_\pi$ and if $x_1,\,x_2 \in S_\pi$ then either $\Delta_{x_1}^\pi = \Delta_{x_2}^\pi$ or
$\Delta_{x_1}^\pi \cap \Delta_{x_2}^\pi = \varnothing$. Put
\[ \mathcal{N}_\pi = \{ N \in \mathcal{F} : \mbox{$\Delta_x^\pi \subset N$ or $\Delta_x^\pi \cap N = \varnothing$
for all $x \in S_\pi$} \}\;; \]
then $\mathcal{N}_\pi$ is a sub-$\sigma$-algebra of $\mathcal{F}$ and $\Delta_x^\pi \in \mathcal{N}_\pi$ for each 
$x \in S_\pi$.

\begin{lemma}\label{lemma_wip_71}
The kernel $\pi$ is $\mathcal{N}_\pi$-measurable, and so in particular $\mathcal{S}_\pi \subset \mathcal{N}_\pi$.
\end{lemma}

\proof
Let $f \in \bMap(\mathcal{F})$, $B \in \Borelpi$ and put $E = \{ x \in X : \pi(f)(x) \in B \}$. Consider 
$x \in E$ and suppose $\Delta_x^\pi \cap E \ne \varnothing$; there thus exists $y \in \Delta_x^\pi \cap E$ and 
so $\varepsilon_y\pi = \varepsilon_x\pi$ and $(\varepsilon_x\pi)(f) = \pi(f)(y) \in B$. Let 
$z \in \Delta_x^\pi$; then $\varepsilon_z\pi = \varepsilon_x\pi = \varepsilon_y\pi$ and hence
$\pi(f)(z) =  (\varepsilon_z\pi)(f) = (\varepsilon_y\pi)(f) \in B$. This shows that $\Delta_x^\pi \subset E$ 
and so $E \in \mathcal{N}_\pi$. It follows that $\pi$ is $\mathcal{N}_\pi$-measurable. \eop

\begin{lemma}\label{lemma_wip_81}
$\pi$ is a proper $\mathcal{N}_\pi$-measurable kernel if and only if $(\varepsilon_x\pi)(\uDelta_x^\pi) = 1$
for each $x \in S_\pi$. 
\end{lemma}

\proof
If $\pi$ is a proper $\mathcal{N}_\pi$-measurable kernel then for each $x \in S_\pi$ we have
$(\varepsilon_x\pi)(\uDelta_x^\pi) = \pi(\uDelta_x^\pi)(x) = \uDelta_x^\pi(x)\pi(1)(x) = 1$ since 
$\uDelta_x^\pi \in \bMap(\mathcal{N}_\pi)$. Thus suppose conversely that
$(\varepsilon_x\pi)(\uDelta_x^\pi) = 1$ for each $x \in S_\pi$. Let $x \in S_\pi$ and $N \in \mathcal{N}_\pi$; 
if $\Delta_x^\pi \subset N$ then $1 = (\varepsilon_x\pi)(\uDelta_x^\pi) \le (\varepsilon_x\pi)(I_N)$ and so 
$(\varepsilon_x\pi)(I_N) = 1 = I_N(x)$, since $x \in \Delta_x^\pi \subset N$. On the other hand, 
$(\varepsilon_x\pi)(I_N) = (\varepsilon_x\pi)(\uDelta_x^\pi I_N) = 0$ if $\Delta_x^\pi \cap N = \varnothing$ 
and hence here $(\varepsilon_x\pi)(I_N) = 0 = I_N(x)$, since $x \notin N$. Therefore in both cases 
$\pi(I_N)(x) = (\varepsilon_x\pi)(I_N) = I_N(x)$, and thus by Lemmas \ref{lemma_wip_61} and 
\ref{lemma_wip_71} $\pi$ is a proper $\mathcal{N}_\pi$-measurable kernel.
\eop

Put $\mathcal{E}_\pi = \mathcal{E} \cap \mathcal{N}_\pi$; thus
$\mathcal{E}_\pi = \{ E \in \mathcal{E} : \mbox{$\Delta_x^\pi \subset E$ or $\Delta_x^\pi \cap E = \varnothing$
for all $x \in S_\pi$} \}$ is a sub-$\sigma$-algebra of $\mathcal{F}$ and by Lemma~\ref{lemma_wip_71}
$\mathcal{S}_\pi \subset \mathcal{E}_\pi \subset \mathcal{E}$. Moreover, $\pi$ is an $\mathcal{E}_\pi$-measurable 
kernel.

\begin{theorem}\label{theorem_wip_31}
The following statements are equivalent:
\begin{evlist}{10pt}{6pt}
\item[(1)]
$\pi$ is a normal $\mathcal{E}$-measurable kernel.

\item[(2)]
$\pi$ is an adapted $\mathcal{E}$-measurable kernel and $(\varepsilon_x\pi)(\uDelta_x^\pi) = 1$ for each 
$x \in S_\pi$. 

\item[(3)]
$\pi$ is an adapted $\mathcal{E}$-measurable kernel and a proper $\mathcal{N}_\pi$-measurable kernel.

\item[(4)]
$\pi$ is a proper $\mathcal{E}_\pi$-measurable kernel and 
$\mathcal{J}_{\mathcal{E}_\pi}(\pi) = \mathcal{J}_{\mathcal{E}}(\pi)$.

\item[(5)]
$\pi$ is a proper $\mathcal{S}_\pi$-measurable kernel and 
$\mathcal{J}_{\mathcal{S}_\pi}(\pi) = \mathcal{J}_{\mathcal{E}}(\pi)$.
\end{evlist}
\end{theorem}

\proof 
(1) $\Leftrightarrow$ (2):\enskip Proposition~\ref{prop_wip_21}.

(2) $\Leftrightarrow$ (3):\enskip Lemma~\ref{lemma_wip_81}.

(3) $\Rightarrow$ (4):\enskip
Since $\pi$ is $\mathcal{E}_\pi$-measurable and $\mathcal{E}_\pi \subset \mathcal{N}_\pi$ it follows that $\pi$
is a proper $\mathcal{E}_\pi$-measurable kernel. Moreover, 
$\mathcal{J}_{\mathcal{E}}(\pi) \subset \mathcal{J}_{\mathcal{E}_\pi}(\pi) \subset \mathcal{J}_\star(\pi)$, since 
$\mathcal{E}_\pi \subset \mathcal{E}$, and by Proposition~\ref{prop_wip_11}~(1) 
$\mathcal{J}_{\mathcal{E}}(\pi) = \mathcal{J}_\star(\pi)$. Hence 
$\mathcal{J}_{\mathcal{E}_\pi}(\pi) = \mathcal{J}_{\mathcal{E}}(\pi)$.

\medskip
(4) $\Rightarrow$ (2):\enskip
Lemma~\ref{lemma_wip_21} implies that $\pi$ is a normal $\mathcal{E}_\pi$-measurable kernel, and thus by 
Proposition~\ref{prop_wip_21} $\pi$ is an adapted $\mathcal{E}_\pi$-measurable kernel and 
$(\varepsilon_x\pi)(\uDelta_x^\pi) = 1$ for all $x \in S_\pi$. But 
$\mathcal{J}_{\mathcal{E}_\pi}(\pi) = \mathcal{J}_{\mathcal{E}}(\pi)$ and so $\pi$ is also an adapted 
$\mathcal{E}$-measurable kernel.

(3) $\Rightarrow$ (5) and (5) $\Rightarrow$ (2) are the same as 
(3) $\Rightarrow$ (4) and (4) $\Rightarrow$ (2), since
$\pi$ is $\mathcal{S}_\pi$-measurable and $\mathcal{S}_\pi \subset \mathcal{N}_\pi$.
\eop

As mentioned near the beginning, there is a further property equivalent to that of being normal, namely that 
$\varepsilon_x\pi \in \ext \mathcal{J}_{\mathcal{E}}(\pi)$ for all $x \in S_\pi$. The equivalence follows 
immediately from the well-known fact that an element of $\mathcal{J}_\mathcal{E}(\pi)$ is extreme if and only if
it is trivial on $\mathcal{E}$. For the sake of completeness we now present a proof of this:

\begin{proposition}\label{prop_wip_31}
An element of $\mathcal{J}_{\mathcal{E}}(\pi)$ is extreme if and only if it is trivial on the $\sigma$-algebra 
$\mathcal{E}$.
\end{proposition}

\proof 
We first need a lemma. Note that if $\mu \in \Prob(\mathcal{F})$ and $h \in \bMap(\mathcal{F})$ with 
$\mu(h) = 1$ then the measure $\measmap{\mu}{h}$ (given by $(\measmap{\mu}{h})(f) = \mu(hf)$ for each 
$f \in \bMap(\mathcal{F})$) is also a probability measure.

\begin{lemma}\label{lemma_wip_91}
Let $\mu \in \mathcal{J}_{\mathcal{E}}(\pi)$ and $h \in \bMap(\mathcal{F})$ with $\mu(h) = 1$. Then 
$\measmap{\mu}{h} \in \mathcal{J}_{\mathcal{E}}(\pi)$ if and only if there exists $h' \in \bMap(\mathcal{E})$
such that $\mu(hf) = \mu(h'f)$ for all $f \in \bMap(\mathcal{F})$. Moreover, in this case we can take 
$h' = \pi(h)$.
\end{lemma}

\proof
If there exists $h' \in \bMap(\mathcal{E})$ with $\mu(hf) = \mu(h'f)$ for all $f \in \bMap(\mathcal{F})$ then
\[
  (\measmap{\mu}{h})(g\pi(f)) =  \mu (hg \pi(f)) = \mu (h'g \pi(f)) 
 = \mu(h'gf) = \mu (hgf) = (\measmap{\mu}{h})(gf)
\]
for all $f \in \bMap(\mathcal{F})$, $g \in \bMap(\mathcal{E})$ and hence 
$\measmap{\mu}{h} \in \mathcal{J}_{\mathcal{E}}(\pi)$. Conversely, if 
$\measmap{\mu}{h} \in \mathcal{J}_{\mathcal{E}}(\pi)$ and $f \in \bMap(\mathcal{F})$ then 
$\mu( h \pi(f)) = \mu( \pi(h) \pi(f)) = \mu(\pi(h)f)$ (since $\mu \in \mathcal{J}_{\mathcal{E}}(\pi)$ and 
$\pi(h),\,\pi(f) \in \bMap(\mathcal{E})$) and therefore
\[ \mu(hf) = (\measmap{\mu}{h})(f) = (\measmap{\mu}{h})(\pi(f)) = \mu( h \pi(f)) = \mu(\pi(h)f)\;.\ \eop \]

If $\mu \in  \mathcal{J}_{\mathcal{E}}(\pi)$ is not extreme then there exist 
$\mu_1,\,\mu_2 \in \mathcal{J}_{\mathcal{E}}(\pi)$ with $\mu_1 \ne \mu_2$ and $0 < a < 1$ such that 
$\mu = a \mu_1 + (1 - a)\mu_2$. Then $\mu_1 \ll \mu$ and so by the Radon-Nikodym Theorem there exists 
$h  \in \bMap(\mathcal{F})$ with $\mu_1 = \measmap{\mu}{h}$, and $\mu(h) = \mu_1(1) = 1$. Therefore by 
Lemma~\ref{lemma_wip_91} there exists $h' \in \bMap(\mathcal{E})$ such that $\mu(h'f) = \mu(hf)$ for all 
$f \in \bMap(\mathcal{F})$, and in particular this implies  $\mu(h') = \mu(h) = 1$. Now $\mu \ne \mu_1$ and so 
let $f \in \bMap(\mathcal{F})$ with  $\mu(f) \ne \mu_1(f) = \mu(hf)$. Then
$\mu(h'f) = \mu(hf) \ne \mu(f) = \mu(h')\mu(f)$ and hence $\mu$ is not trivial on $\mathcal{E}$.

Conversely, suppose $\mu$ is not trivial on $\mathcal{E}$, and therefore there exists $E \in \mathcal{E}$ with 
$0 < \mu(I_E) < 1$. Put $a = \mu(I_E)$ and let $\mu_1 = \measmap{\mu}{g_1}$ and  $\mu_2 = \measmap{\mu}{g_2}$ 
with $g_1 = a^{-1}I_E$ and $g_2 = (1 - a)^{-1}I_{X\setminus E}$. Then $g_j \in \bMap(\mathcal{E})$ and
$\mu(g_j) = 1$ and so by Lemma~\ref{lemma_wip_91} $\mu_j \in \mathcal{J}_{\mathcal{E}}(\pi)$ for $j = 1,\,2$.
However $\mu = a\mu_1 + (1- a)\mu_2$ and clearly $\mu_1 \ne \mu_2$; hence 
$\mu \notin \ext \mathcal{J}_{\mathcal{E}}(\pi)$. This completes the proof of Proposition~\ref{prop_wip_31}.
\eop

We next present a condition which is equivalent to the existence of a proper $\mathcal{E}$-measurable refinement.

If $\mathcal{D}$ is a sub-$\sigma$-algebra of $\mathcal{F}$ and $A \subset X$ then the
\definition{trace $\sigma$-algebra $\mathcal{D}_{|A}$} is the $\sigma$-algebra of subsets of $A$ given by
$\{ B \subset A: \mbox{$B = A \cap D$ for some $D \in \mathcal{D}$} \}$. If $A \in \mathcal{D}$ then 
$\mathcal{D}_{|A}$ just consists of the subsets of $A$ which lies in $\mathcal{D}$. If $\mathcal{D}$ is countably 
generated then the trace $\sigma$-algebra $\mathcal{D}_{|A}$ is also countably generated for each $A \subset X$.
The next fact is essentially part of Theorem 1 in Blackwell and Dubins \cite{blackwelldubins}:

\begin{lemma}\label{lemma_wip_101}
If $\pi$ is proper then the trace $\sigma$-algebra $\mathcal{E}_{|S_\pi}$ is countably generated.
\end{lemma}

\proof 
If $E \in \mathcal{E}$ then $I_{E \cap S_\pi} = I_E\pi(1) = \pi(I_E)  \in  \bMap(\mathcal{S}_\pi)$ and hence 
$E \cap S_\pi \in \mathcal{S}_\pi$. Thus $E \cap S_\pi = (E \cap S_\pi) \cap S_\pi \in (\mathcal{S}_\pi)_{|S_\pi}$ 
for all $E \in \mathcal{E}$, which implies  $\mathcal{E}_{|S_\pi} \subset (\mathcal{S}_\pi)_{|S_\pi}$.
Moreover  $(\mathcal{S}_\pi)_{|S_\pi} \subset \mathcal{E}_{|S_\pi}$, since $\mathcal{S}_\pi \subset \mathcal{E}$,
and hence $\mathcal{E}_{|S_\pi} = (\mathcal{S}_\pi)_{|S_\pi}$. But $(\mathcal{S}_\pi)_{|S_\pi}$ is countably 
generated, since $\mathcal{S}_\pi$ is, and therefore $\mathcal{E}_{|S_\pi}$ is countably generated. \eop

Consider the special case in which $\pi$ is a probability kernel (and so $S_\pi = X$). Then in particular 
Lemma~\ref{lemma_wip_101} shows that $\mathcal{E}$ is countably generated. This implies that if $\mathcal{E}$ 
is not countably generated then there are no proper $\mathcal{E}$-measurable probability kernels.

Let us say that a set $D$ is \definition{$\pi$-full} if $D\in \mathcal{E}$ and $\mu(I_D) = 1$ for all 
$\mu \in \mathcal{J}_{\mathcal{E}}(\pi)$

\begin{theorem}\label{theorem_wip_41}
The following statements are equivalent:
\begin{evlist}{10pt}{6pt}
\item[(1)]
There exists a proper $\mathcal{E}$-measurable refinement of $\pi$.

\item[(2)]
The trace $\sigma$-algebra $\mathcal{E}_{|D}$ is countably generated for some $\pi$-full set $D$.
\end{evlist}
\end{theorem}

\proof 
(2) $\Rightarrow$ (1):\enskip
Since $D \cap S_\pi$ is also $\pi$-full and $\mathcal{E}_{|D \cap S_\pi}$ is countably generated, we can assume 
that $D \subset S_\pi$. Extend each mapping $g \in \bMap(D)$ to a mapping $g^* \in \bMap(X)$ by putting 
$g^*(x) = 0$ for all $x \in X \setminus D$; thus if $g \in \bMap(\mathcal{E}_{|D})$ then 
$g^* \in \bMap(\mathcal{E})$, since $D \in \mathcal{E}$. As well as the countable determining set $G$ for  
finite measures on $(X,\mathcal{F})$  choose a countable determining set $G'$ for finite measures on 
$(X,\mathcal{E}_{|D})$. Let
\[ C = \{ x \in D : \pi(g^*f)(x) = g^*(x)\pi(f)(x)\ \mbox{for all $g \in G'$, $f \in G$} \}\]
then $C = \bigcap_{f \in G} \bigcap_{g \in G'} C_{f,g}$, where 
$C_{f,g} =  \{ x \in D: \pi(g^*f)(x) = g^*(x)\pi(f)(x) \}$. In particular $C \in \mathcal{E}$. Let 
$\mu \in \mathcal{J}_{\mathcal{E}}(\pi)$, $g \in G'$, $f \in G$; then for all $h \in \bMap(\mathcal{E})$ we have
$\mu(h\pi(g^*f)) = \mu(hg^*f)) = \mu(hg^*\pi(f))$, and so $\mu(I_{C_{g,h}}) = 1$, since $\mu(I_D) = 1$. Hence 
$\mu(I_C) = 1$, since $G \times G'$ is countable, which shows that $C$ is $\pi$-full.

Now let $\varrho$ be the restriction of $\pi$ to $C$; thus by Lemma~\ref{lemma_wip_11} $\varrho$ is a refinement 
of $\pi$. Exactly as in the proof of Theorem~\ref{theorem_wip_11} it follows that
$\varrho(g^*f) = g^*\varrho(f)$ for all $f \in \bMap(\mathcal{F})$ and $g \in \bMap(\mathcal{E}_{|D})$. In 
particular, if $x \in C = S_\varrho$ then we have 
$(\varepsilon_x\varrho)(I_D) = \varrho(I_D)(x) = I_D(x)\varrho(1)(x) = 1$, since $I_D= I_D^*$ and 
$S_\varrho \subset D$. Let $g \in \bMap(\mathcal{E})$ and $h$ be the restriction of $g$ to $S_D$; then 
$h \in \bMap(\mathcal{E}_{|D})$ and $h^* = I_Dg$. Thus if $x \in C$ then for all $f \in \bMap(\mathcal{F})$
\begin{eqnarray*}
\varrho(gf)(x) = (\varepsilon_x\varrho)(gf) &=& (\varepsilon_x\varrho)(I_{D}gf)\\ 
&=&(\varepsilon_x\varrho)(h^*f) 
= \varrho(h^*f)(x) = h^*(x)\varrho(f)(x) = g(x)\varrho(f)(x)\;.
\end{eqnarray*}
On the other hand, if $x \in X \setminus C$ then $\varrho(gf)(x) = 0 = g(x)\varrho(f)(x)$, which shows that 
$\varrho(gf) = g\varrho(f)$ for all $g \in \bMap(\mathcal{E})$, $f \in \bMap(\mathcal{F})$ and therefore that 
$\varrho$ is a proper $\mathcal{E}$-measurable kernel.

\medskip
(1) $\Rightarrow$ (2):\enskip
If $\varrho$ is a proper $\mathcal{E}$-measurable refinement of $\pi$ then by Lemma~\ref{lemma_wip_101} the 
trace $\sigma$-algebra $\mathcal{E}_{|S_\varrho}$ is countably generated. Moreover $S_\varrho$ is $\varrho$-full 
and thus $\pi$-full (since $\mathcal{J}_{\mathcal{E}}(\pi) = \mathcal{J}_{\mathcal{E}}(\varrho)$). \eop

We now apply Theorems \ref{theorem_wip_11} and \ref{theorem_wip_21} to present a prototype of a result which 
occurs in Dynkin's construction of an entrance boundary (Dynkin \cite{dynkin1} and \cite{dynkin2}) and in 
F\"ollmer's representation of Gibbs states on the tail $\sigma$-algebra (F\"ollmer \cite{foellmer}). The 
result given below (Theorem~\ref{theorem_wip_51}) is a modification of the account to be found in
Preston \cite{preston}.

In what follows let $\{\mathcal{E}_n\}_{n\ge 1}$ be a decreasing sequence of sub-$\sigma$-algebras of 
$\mathcal{F}$ and for each $n \ge 1$ let $\pi_n$ be an $\mathcal{E}_n$-measurable quasi-probability kernel
such that the sequence $\{\mathcal{J}_{\mathcal{E}_n}(\pi_n)\}_{n\ge 1}$ is decreasing, i.e.,
$\mathcal{J}_{\mathcal{E}_n}(\pi_n) \subset \mathcal{J}_{\mathcal{E}_m}(\pi_m)$ whenever $m \le n$. 

We suppose for each $n \ge 1$ that either the $\sigma$-algebra $\mathcal{E}_n$ is countably generated or that
$\mathcal{J}_{\mathcal{E}}(\pi_n) = \mathcal{J}_\star(\pi_n)$. In the first case there exists by 
Theorem~\ref{theorem_wip_11} a proper $\mathcal{E}_n$-measurable refinement $\varrho_n$ of $\pi_n$. In the 
second Theorem~\ref{theorem_wip_21} guarantees the existence of a normal refinement $\varrho_n$ of $\pi_n$.
Since $\mathcal{J}_{\mathcal{E}}(\varrho_n) = \mathcal{J}_{\mathcal{E}}(\pi_n)$ for each $n$ the sequence 
$\{\mathcal{J}_{\mathcal{E}_n}(\varrho_n)\}_{n\ge 1}$ is also decreasing. Moreover, $\varrho_n$ is normal for each 
$n$ and therefore by Proposition~\ref{prop_wip_11}~(2) the sequence of kernels $\{\varrho_n\}_{n\ge 1}$ is 
compatible in that $\varrho_n\varrho_m = \varrho_n$ whenever $m \le n$. Let 
$\mathcal{E} = \bigcap_{n \ge 1} \mathcal{E}_n$.

\begin{theorem}\label{theorem_wip_51}
Suppose that $(X,\mathcal{F})$ is a standard Borel space. Then there exists a normal 
$\mathcal{E}$-measurable quasi-probability kernel $\varrho$ such that
\[ \bigcap_{n \ge 1} \mathcal{J}_{\mathcal{E}_n}(\pi_n) = \mathcal{J}_{\mathcal{E}}(\varrho)\;.\]
\end{theorem}

\proof
This is divided into two parts. The first only requires $\mathcal{F}$ to be countably generated and shows that 
if there exists an $\mathcal{E}$-measurable quasi-probability kernel $\pi$ such that
$\bigcap_{n \ge 1} \mathcal{J}_{\mathcal{E}_n}(\pi_n) \subset \mathcal{J}_\mathcal{E}(\pi)$ then there exists a 
restriction $\varrho$ of $\pi$ which is normal and such that 
$\bigcap_{n \ge 1} \mathcal{J}_{\mathcal{E}_n}(\pi_n) = \mathcal{J}_{\mathcal{E}}(\varrho)$. The second part shows 
that if $(X,\mathcal{F})$ is a standard Borel space then there does exist an 
$\mathcal{E}$-measurable quasi-probability kernel $\pi$ such that
$\bigcap_{n \ge 1} \mathcal{J}_{\mathcal{E}_n}(\pi_n) \subset \mathcal{J}_\mathcal{E}(\pi)$.

\begin{proposition}\label{prop_wip_41}
If $\pi$ is an $\mathcal{E}$-measurable quasi-probability kernel such that
\[ \bigcap_{n \ge 1} \mathcal{J}_{\mathcal{E}_n}(\pi_n) \subset \mathcal{J}_\mathcal{E}(\pi)\]
then there exists a restriction $\varrho$ of $\pi$ which is normal and such that
\[ \bigcap_{n \ge 1} \mathcal{J}_{\mathcal{E}_n}(\pi_n) = \mathcal{J}_{\mathcal{E}}(\varrho)\;.\]
\end{proposition}

\proof
For each $n \ge 1$ let $\varrho_n$ be the normal refinement of $\pi_n$ introduced above. Then 
$\mathcal{J}_{\mathcal{E}_n}(\varrho_n) = \mathcal{J}_{\mathcal{E}_n}(\pi_n)$ and by 
Proposition~\ref{prop_wip_11}~(2) $\mathcal{J}_{\mathcal{E}_n}(\varrho_n) = \mathcal{J}_\star(\varrho_n)$.
Without loss of generality we can thus assume that 
$\mathcal{J}_{\mathcal{E}}(\pi_n) = \mathcal{J}_\star(\pi_n)$ for each $n \ge 1$. Put 
$\mathcal{J} = \bigcap_{n \ge 1} \mathcal{J}_{\mathcal{E}_n}(\pi_n)$ and let 
$D = \{ x \in S_\pi : \varepsilon_x\pi \in \mathcal{J}\}$; thus $D = \bigcap_{n \ge 1} D_n$, where 
\begin{eqnarray*}
D_n &=& \{ x \in S_\pi : \varepsilon_x\pi \in \mathcal{J}_{\mathcal{E}_n}(\pi_n)\}\\
                          &=& \{ x \in S_\pi : \varepsilon_x\pi \in \mathcal{J}_\star(\pi_n)\}
 = \{ x \in S_\pi : \varepsilon_x\pi = \varepsilon_x(\pi\pi_n) \}
\end{eqnarray*}
and it follows that $D_n = \bigcap_{f \in G} D^f_n$, where 
\[ D^f_n  =\{ x \in S_\pi : (\varepsilon_x\pi)(f) = (\varepsilon_x(\pi\pi_n))(f) \}
   = \{ x \in S_\pi : \pi(f)(x) = \pi(\pi_n(f))(x) \} \;. \]
In particular this implies that $D \in \mathcal{E}$. We first show that $\mu(I_D) = 1$ for all 
$\mu \in \mathcal{J}$. If $\mu \in \mathcal{J}_{\mathcal{E}}(\pi) \cap \mathcal{J}_{\mathcal{E}_n}(\pi_n)$ then 
$\mu(g\pi(\pi_n(f)) = \mu(g \pi_n(f)) = \mu(gf)  = \mu(g \pi(f))$ for all $g \in \bMap(\mathcal{E})$, and hence 
$\mu(I_{D^f_n}) = 1$, since $\mu(I_{S_\pi}) = 1$. Therefore $\mu(I_{D_n}) = 1$, since $G$ is countable, i.e., 
$\mu(I_{D_n}) = 1$ for all $\mu \in \mathcal{J}_{\mathcal{E}}(\pi) \cap \mathcal{J}_{\mathcal{E}_n}(\pi_n)$. But 
$D = \bigcap_{n\ge 1} D_n$ and 
$\mathcal{J} = \bigcap_{n \ge 1} \bigl(\mathcal{J}_{\mathcal{E}}(\pi)\cap \mathcal{J}_{\mathcal{E}_n}(\pi_n)\bigr)$,
and so $\mu(I_D) = 1$ for all $\mu \in \mathcal{J}$.

Let $\tau$ be the restriction of $\pi$ to $D$; we next show that $\mathcal{J}_\star(\tau) \subset \mathcal{J}$.
Let $n \ge 1$,  $f \in \bMap(\mathcal{F})$ and $g \in \bMap(\mathcal{E}_n)$; if $x \in S_{\tau} = D$ then 
$\varepsilon_x\tau = \varepsilon_x\pi \in  \mathcal{J}_{\mathcal{E}_n}(\pi_n)$ and so 
$(\varepsilon_x\tau)(g\pi_n(f)) = (\varepsilon_x\tau)(gf)$. But this also holds trivially when 
$x \in X \setminus D$, since then $\varepsilon_x\tau = 0$, and hence $\tau(g\pi_n(f)) = \tau(gf)$. Thus if 
$\mu \in \mathcal{J}_\star(\tau)$ then
\[ \mu(g\pi_n(f)) = (\mu\tau)(g\pi_n(f))= \mu(\tau(g\pi_n(f))) = \mu(\tau(gf)) = (\mu\tau)(gf) = \mu(gf)\;, \]
which implies that $\mu \in \mathcal{J}_{\mathcal{E}_n}(\pi_n)$. Therefore
$\mathcal{J}_\star(\mathcal{\tau}) \subset \mathcal{J}_{\mathcal{E}_n}(\pi_n)$ for each $n \ge 1$ and so 
$\mathcal{J}_\star(\mathcal{\tau}) \subset \mathcal{J}$. 

Now if $\mu \in \mathcal{J}$ then $\mu \in \mathcal{J}_{\mathcal{E}}(\pi)$ and it was shown above that 
$\mu(I_D) = 1$; therefore $\mu(gf) = \mu(g\pi(f))= \mu(gI_D\pi(f)) = \mu(g\tau(f))$ for all  
$g \in \bMap(\mathcal{E})$, $f \in \bMap(\mathcal{F})$ and hence $\mu \in \mathcal{J}_{\mathcal{E}}(\tau)$. 
This shows that $\mathcal{J} \subset \mathcal{J}_{\mathcal{E}}(\tau)$. 

Combining the above inclusions gives
$\mathcal{J}_\star(\mathcal{\tau}) \subset \mathcal{J} \subset \mathcal{J}_{\mathcal{E}}(\tau)
       \subset \mathcal{J}_\star(\mathcal{\tau})$, 
which means that $\mathcal{J} = \mathcal{J}_{\mathcal{E}}(\tau) = \mathcal{J}_\star(\mathcal{\tau})$. Finally, 
by Theorem~\ref{theorem_wip_21} there exists a normal refinement $\varrho$ of $\tau$, since 
$\mathcal{J}_{\mathcal{E}}(\tau) = \mathcal{J}_\star(\mathcal{\tau})$, and $\varrho$ is the required restriction 
of $\pi$. \eop

\begin{proposition}\label{prop_wip_51}
Suppose $(X,\mathcal{F})$ is a standard Borel space. Then there exists an $\mathcal{E}$-measurable 
quasi-probability kernel $\varrho$ such that
\[ \bigcap_{n \ge 1} \mathcal{J}_{\mathcal{E}_n}(\pi_n) \subset \mathcal{J}_{\mathcal{E}}(\varrho)\;.\]
\end{proposition}

\proof
Put $\mathcal{J} = \bigcap_{n \ge 1} \mathcal{J}_{\mathcal{E}_n}(\pi_n)$. Since $(X,\mathcal{F})$ is a standard 
Borel space there exists what Dynkin \cite{dynkin1} calls a support system. This is a countable determining set 
$G$ for finite measures on $(X,\mathcal{F})$ having the additional property that if $\{\mu_n\}_{n\ge 1}$ is a 
sequence from $\Prob(\mathcal{F})$ such that $\lim_n \mu_n(f)$ exists for each $f \in G$ then there exists 
(a unique) $\mu \in \Prob(\mathcal{F})$ such that $\lim_n \mu_n(f) = \mu(f)$ for all $f \in G$. Let
\[X_0 = \bigl\{ x \in X: \lim\limits_{n\to\infty} \pi_n(f)(x)\ \mbox{exists for all}\ f \in G \bigr\}\;,\]
We first show that $X_0 \in \mathcal{E}$ and $\mu(I_{X_0}) = 1$ for each $\mu \in \mathcal{J}$. For each 
$f \in G$ let $X_f$ denote the set of those elements $x \in X$ for which the limit $\lim_{n} \pi_n(f)(x)$ exists, 
thus $X_0 = \bigcap_{f \in G} X_f$ and therefore, since $G$ is countable, it is enough to show for each 
$f \in G$ that $X_f \in \mathcal{E}$ and $\mu(I_{X_f}) = 1$ for each $\mu \in \mathcal{J}$. Now since $X_f$ 
doesn't depend on the first $m$ terms of the sequence $\{\pi_n(f)\}_{n\ge 1}$ for any $m \ge 1$ it follows that
$X_f \in \mathcal{E}$. Let $\mu \in \mathcal{J}$; then $\mu(h \pi_n(f)) = \mu(hf)$ for all 
$h \in \bMap(\mathcal{E}_n)$, which means that $\pi_n(f)$ is a version of the conditional expectation of $f$ with 
respect to $\mathcal{E}_n$ (with the measure $\mu$ here fixed). Thus by the backward martingale convergence 
theorem $\mu(I_{X_f}) = 1$. For each $f \in G$ there is an element $\tau_f \in \bMap(\mathcal{E})$ given by
\[ \tau_f(x) =  \left\{ \begin{array}{cl}
                     \lim\limits_{n\to \infty} \pi_n(f)(x) &\ \mbox{if}\ x \in X_0\;,\\
                                 0                 &\ \mbox{if}\ x \in X \setminus X_0\;.
                \end{array} \right. \]
                            
Note that if $x \in X_0$ then $\tau_1(x) \in \{0,1\}$, since $\pi_n(1)(x) \in \{0,1\}$ for all $n$. Let 
$X_1 = \{ x \in X_0 : \tau_1(x) = 1 \}$, so $X_1 \in \mathcal{E}$. If $x \in X_1$ then
$(\varepsilon_x\pi_n)(1) = \pi_n(1)(x) = 1$ for all $n \ge n_x$ for some $n_x \ge 1$, which means that
$\{\varepsilon_x\pi_n\}_{n \ge n_x}$ is a sequence of probability measures with 
$\lim_n (\varepsilon_x\pi_n)(f) = \tau_f(x)$ for all $f \in G$. Thus, since $G$ is a support system, there exists 
$\mu_x \in \Prob(\mathcal{F})$ such that $\mu_x(f) = \tau_f(x)$ for all $f \in G$. Now define 
$\pi :\bMap(\mathcal{F}) \to \bMap(X)$ by
\[ \pi(f)(x) =  \left\{ \begin{array}{cl}
                                         \mu_x(f)  &\ \mbox{if}\ x \in X_1\;,\\
                                 0                 &\ \mbox{if}\ x \in X \setminus X_1\;.
                \end{array} \right. \]
Then $\pi(f) = \tau_f \in \bMap(\mathcal{E})$ for all $f \in G$ (noting that
$\pi(f)(x) = 0 = \tau_f(x)$ for all $x \in X_0 \setminus X_1$) and it is straightforward to show that $\pi$ is 
an $\mathcal{E}$-measurable quasi-probability kernel. Let $\mu \in \mathcal{J}$; if $g \in \bMap(\mathcal{E})$ 
and $f \in G$ then $\mu(gf) = \mu(g \pi_n(f))$ for each $n \ge 1$ and $\mu(I_{X_0}) = 1$; thus by the dominated 
convergence theorem $\mu(g f) = \lim_{n} \mu(g I_{X_0}\pi_n(f)) = \mu(g \tau_f) = \mu(g \pi(f))$. Since $G$ is a 
generator for $\bMap(\mathcal{F})$ it follows that $\mu(g f) = \mu(g \pi(f))$ for all $f \in \bMap(\mathcal{F})$, 
which implies that $\mu \in \mathcal{J}_{\mathcal{E}}(\pi)$. Therefore 
$\mathcal{J} \subset \mathcal{J}_{\mathcal{E}}(\pi)$. \eop

Theorem~\ref{theorem_wip_51} now follows immediately from Propositions \ref{prop_wip_41} and \ref{prop_wip_51}. 
\eop

As a final topic we look at a problem considered in Sokal \cite{sokal}, Goldstein \cite{goldstein} and 
Preston \cite{preston2}: Again let $\{\mathcal{E}_n\}_{n\ge 1}$ be a decreasing sequence of sub-$\sigma$-algebras 
of $\mathcal{F}$ and assume now that $\mathcal{E}_n$ is countably generated for each $n \ge 1$. Fix 
$\mu \in \Prob(\mathcal{F})$. Does there then exist a proper $\mathcal{E}_n$-measurable quasi-probability kernel
$\varrho_n$ for each $n \ge 1$ such that $\varrho_n\varrho_m = \varrho_n$ whenever $m \le n$ and
$\mu \in \mathcal{J}_{\mathcal{E}_n}(\varrho_n)$ for all $n \ge 1$? We will see that if $(X,\mathcal{F})$ is a 
standard Borel space then there do exist kernels with these properties; this follows more-or-less directly from 
the arguments  found in Sokal \cite{sokal}.

\begin{lemma}\label{lemma_wip_111}
Let $\mathcal{E}$ and $\mathcal{E}'$ be sub-$\sigma$-algebras of $\mathcal{F}$ with 
$\mathcal{E} \subset \mathcal{E}'$ and $\mathcal{E}$ countably generated. Suppose  
$\mu \in \mathcal{J}_{\mathcal{E}}(\pi) \cap \mathcal{J}_{\mathcal{E}'}(\pi')$, where $\pi$ is an 
$\mathcal{E}$-measurable and $\pi'$ an $\mathcal{E}'$-measurable quasi-probability kernel. Then there exists a 
proper $\mathcal{E}$-measurable quasi-probability kernel $\varrho$ with $\varrho\pi' = \varrho$ and 
$\mu \in \mathcal{J}_{\mathcal{E}}(\varrho)$.
\end{lemma}

\proof
Let $D = \{ x \in X : (\pi\pi')(f)(x) = \pi(f)(x)\ \mbox{for all $f \in \bMap(\mathcal{F})$} \}$;  by 
Lemma~\ref{lemma_wip_31} $D = \bigcap_{f \in G} D_f$, where $D_f = \{ x \in X : \pi(\pi'(f))(x) = \pi(f)(x) \}$ 
and so in particular $D \in \mathcal{E}$. If $g \in \bMap(\mathcal{E})$  then
$\mu(g\pi(\pi'(f))) = \mu(g\pi'(f)) = \mu(gf) = \mu(g(\pi(f))$, since 
$\mu \in \mathcal{J}_{\mathcal{E}}(\pi) \cap \mathcal{J}_{\mathcal{E}'}(\pi')$, and thus
$\mu(g\pi(\pi'(f))) = \mu(g\pi'(f)) = \mu(gf) = \mu(g(\pi(f))$. Hence $\mu(I_{D_f}) = 1$ for each $f \in G$ and 
therefore $\mu(I_D) = 1$. Let $\tau$ be the restriction of $\pi$ to $D$; if $x \in S_\tau = D$ then
$(\tau\pi')(f)(x) = (\pi\pi')(f)(x) = \pi(f)(x) = \tau(f)(x)$ for all $f \in \bMap(\mathcal{F})$ and if 
$x \notin S_\tau$ then $(\tau\pi')(f)(x) = 0 = \tau(f)(x)$, which implies that $\tau\pi' = \tau$. Moreover, 
$\mu(g\tau(f)) = \mu(gI_D\pi(f)) = \mu(g\pi(f)) = \mu(gf)$ for all $f \in \bMap(\mathcal{F})$ and 
$g \in \bMap(\mathcal{E})$, since $\mu(I_D) = 1$, and so $\mu \in \mathcal{J}_{\mathcal{E}}(\tau)$. Now by 
Theorem~\ref{theorem_wip_11} there exists a proper refinement $\varrho$ of $\tau$ and in particular 
$\mu \in \mathcal{J}_{\mathcal{E}}(\tau) = \mathcal{J}_{\mathcal{E}}(\varrho)$. Also by 
Proposition~\ref{prop_wip_11} $\varrho\tau = \varrho$ and hence
$\varrho = \varrho\tau = \varrho(\tau \pi') = (\varrho\tau) \pi' = \varrho\pi'$.
\eop

\begin{theorem}\label{theorem_wip_61}
Suppose that $(X,\mathcal{F})$ is standard Borel. Then for each $n \ge 1$ there exists a proper 
$\mathcal{E}_n$-measurable quasi-probability kernel $\varrho_n$ such that $\varrho_n\varrho_m = \varrho_n$ 
whenever $m \le n$ and $\mu \in \mathcal{J}_{\mathcal{E}_n}(\varrho_n)$ for all $n \ge 1$.
\end{theorem}

\proof 
Since $(X,\mathcal{F})$ is a standard Borel space there exists for each $n \ge 1$ an 
$\mathcal{E}_n$-measurable quasi-probability kernel $\pi_n$ such that 
$\mu \in \mathcal{J}_{\mathcal{E}_n}(\pi_n)$ (Ji\v{r}ina \cite{jirina}). By Theorem~\ref{theorem_wip_11} there 
exists a proper refinement $\varrho_1$ of $\pi_1$ and in particular 
$\mu \in \mathcal{J}_{\mathcal{E}_1}(\varrho_1)$. Let $n \ge 1$ and suppose for $j = 1,\,\ldots,\,n$ there exists 
a proper $\mathcal{E}_j$-measurable quasi-probability kernel $\varrho_j$ such that 
$\mu \in \mathcal{J}_{\mathcal{E}_j}(\varrho_j)$ for all $1 \le j \le n$ and $\varrho_j\varrho_k = \varrho_j$ 
whenever $k \le j \le n$. By Lemma~\ref{lemma_wip_111} there then exists a proper $\mathcal{E}_{n+1}$-measurable 
quasi-probability kernel $\varrho_{n+1}$ such that $\mu \in \mathcal{J}_{\mathcal{E}_{n+1}}(\varrho_{n+1})$ and 
$\varrho_{n+1}\varrho_n = \varrho_{n+1}$. Now if $m \le n$ then 
$\varrho_{n+1} = \varrho_{n+1}\varrho_n = \varrho_{n+1}(\varrho_n\varrho_m)
                         = (\varrho_{n+1}\varrho_n)\varrho_m = \varrho_{n+1}\varrho_m$ 
which implies that $\varrho_j\varrho_k = \varrho_j$ whenever $k \le j \le n+1$, since by 
Proposition~\ref{prop_wip_11} $\varrho_{n+1}\varrho_{n+1} = \varrho_{n+1}$. The result therefore follows by 
induction on $n$. \eop

\bigskip
\bigskip

{\sc Fakult\"at f\"ur Mathematik, Universit\"at Bielefeld}\\
{\sc Postfach 100131, 33501 Bielefeld, Germany}\\
\textit{E-mail address:} \texttt{preston@math.uni-bielefeld.de}\\


\begin{thebibliography}{99}


\bibitem{blackwelldubins} Blackwell, D., Dubins, L. (1975): On existence and non-existence of proper, regular,
conditional distributions. Annals of Probability 3, 741-752.


\bibitem{dynkin1} Dynkin, E. (1971): The initial and final behavior of trajectories of Markov processes.
Russian Math.\ Surveys 26, 165-185.

\bibitem{dynkin2} Dynkin, E. (1978): Sufficient statistics and extreme points. Annals of Probability 6,
705-730.

\bibitem{foellmer} F\"ollmer, H. (1975): Phase transition and Martin boundary. In:
S\'em.\ Prob.\ Strasbourg, SLNM 465.

\bibitem{georgii} Georgii, H.-O. (1988): Gibbs Measures and Phase Transitions. de Gruyter Studies in Mathematics,
9.

\bibitem{goldstein} Goldstein, S. (1978): A note on specifications. ZWvG, 46, 45-51.

\bibitem{halmos} Halmos, P. (1941): The decomposition of measures. Duke Math.\ J., 8, 386-392.

\bibitem{jirina} Ji\v{r}ina, M. (1954): Conditional probabilities on $\sigma$-algebras with countable basis. 
Czech.\ Math.\ J., 79, 372-380.

\bibitem{preston} Preston, C. (1976): Random Fields. SLNM 534.

\bibitem{preston2} Preston, C. (1980): Construction of specification. In: Quantum fields -- algebras,
processes. Wien: Springer.


\bibitem{sokal} Sokal, A. (1981): Existence of compatible families of proper regular conditional probabilities.
  ZWvG 56, 537-548.




\end{thebibliography}
\end{document}